\documentclass [a4paper,11pt]{amsart}
\pdfoutput=1
\usepackage{graphicx}
\usepackage[applemac]{inputenc}
\usepackage[T1]{fontenc}
\usepackage{amssymb}
\usepackage{mathrsfs}
\usepackage{amsfonts,amscd,amsmath}
\usepackage{amssymb,latexsym}
\usepackage[all]{xy}

\newtheorem{thm}{Theorem}[section]
\newtheorem{cor}[thm]{Corollary}
\newtheorem{lem}[thm]{Lemma}
\newtheorem{prop}[thm]{Proposition}

\theoremstyle{definition}
\newtheorem{definition}[thm]{Definition}
\newtheorem{defprop}[thm]{Proposition-definition}
\newtheorem{open}{Open Problem}
\newtheorem{rem}[thm]{Remark}
\newtheorem{ex}[thm]{Example}
\usepackage{mystyle}

\usepackage{color}

\usepackage{hyperref}
\hypersetup{colorlinks,%
citecolor=black,%
filecolor=black,%
linkcolor=black,%
urlcolor=black,%
pdftex}

\begin{document}

\author{Juan J. Nuño-Ballesteros, Guillermo Peñafort-Sanchis}

\address{Departament de Geometria i Topologia,
Universitat de Val\`encia, Campus de Burjassot, 46100 Burjassot,
Spain}

\email{Juan.Nuno@uv.es}

\email{Guillermo.Penafort@uv.es}

\title[Multiple point spaces]{Multiple point spaces\\ of finite holomorphic maps}

\thanks{Work partially supported by DGICYT Grant MTM2015--64013--P}

\begin{abstract}
Given a finite holomorphic map $f:X\to Y$ between manifolds, we show that there exists a unique possible definition of $k$th multiple point space $D^k(f)$ with the following properties: $D^k(f)$ is a closed subspace of $X^k$, $D^k(f)$ is the closure of the set of strict $k$-multiple points when $f$ is stable and $D^k(f)$ is well behaved under deformations. Our construction coincides with Mond's double point space and Mond's $k$th multiple point space for corank one singularities. We also give some interesting properties of the double point space and prove that in many cases it can be computed as the zero locus of the quotient of ideals $(f\times f)^*\mathscr I_{\Delta(Y)}:\mathscr I_{\Delta(X)}$, where $\mathscr I_{\Delta(Z)}$ is the defining ideal sheaf of the diagonal in $Z\times Z$.
\end{abstract}

								\maketitle
													\section*{Introduction}
This work is about the analytic structure of multiple point spaces of finite holomorphic maps $f\colon X\to Y$ between manifolds. A natural definition of the $k$th multiple point space of $f$ is just the closure of the set of strict $k$-multiple points, i.e., 
 the $k$-tuples $(x^{(1)},\dots, x^{(k)})\in X^k$ such that $x^{(i)}\neq x^{(j)}$ and $f(x^{(i)})=f(x^{(j)})$, $ \forall i\ne j$. Unfortunately, this {construction} does not prove useful in the context of singularities of maps, because {it does not behave well under deformations} (see Example \ref{exDeformationOfCusp}). Several authors have worked on multiple point spaces and have introduced different solutions to deal with {this problem} (see for instance \cite{KleimanMultiplePointFormulasI, MararMondCorank1, MondSomeRemarks,  MondPellikaanFittingIdeals, Ronga1972La-classe-duale}). However, the relations between these spaces are not always clear.

In \cite{MondSomeRemarks}, Mond gives a definition of double point space $D^2(f)$ for map germs $f:(\C^n,0)\to(\C^p,0)$, with $n\leq p$, which may be non reduced in general, and whose underlying set is given by the pairs $(x,x')$ such that either $(x,x')$ is a strict double point or $x=x'$ is a singular point. In the same paper, Mond also defines the higher multiple point spaces $D^k(f)$, but only for corank 1 map germs. These spaces are studied in greater depth by Marar-Mond in \cite{MararMondCorank1}. The main result of that paper is that $f$ is stable if and only if all the spaces $D^k(f)$ are smooth of the expected dimension or empty, and that $f$ is finitely determined if and only if all the spaces $D^k(f)$ are isolated complete intersection singularities of the expected dimension or $D^k(f)\subset\{0\}$.

Our approach is based on the ideas of Gaffney in \cite{GaffneyMultiplePointsAndAssociated}. We show that given a finite holomorphic map $f:X\to Y$ between manifolds, there exists only one way to define $D^k(f)$ with the following two properties: first, $D^k(f)$ is the closure of the set of strict $k$-multiple points when $f$ is stable, and second, $D^k(f)$ is well behaved under deformations. This means that if $F$ is an unfolding of $f$, then $D^k(F)$ is a family whose special fibre is $D^k(f)$. Then, we prove that our construction coincides with Mond's definition when $k=2$ or when $f$ has only corank 1 singularities.

The last two sections are dedicated to the double point space $D^2(f)$. We look at interesting geometrical properties of this space and study in which conditions it is smooth, Cohen-Macaulay, reduced or normal. Finally, let $\dim X=n$ and $\dim Y=p$, with $n\leq p$. We prove that if $\dim D^2(f)=2n-p$, and the set of points of corank $\ge 2$ has dimension $<2n-p$, then $D^2(f)$ can be computed as
$$D^2(f)=V((f\times f)^*\mathscr I_{\Delta(Y)}:\mathscr I_{\Delta(X)}),$$
 where $\mathscr I_{\Delta(Z)}$ is the defining ideal sheaf of the diagonal in $Z\times Z$. The result implies, for instance, that the equality holds in the following cases:
\begin{enumerate}
\item $p=n$,
\item $p=n+1$ and $f$ generically one-to-one,
\item $f$ stable,
\item $p<2n$ and $f$ finitely determined (for map germs).
\end{enumerate}

\medskip
\emph{Acknowledgements:}
The authors thank W.L. Marar and D. Mond for their valuable comments and suggestions.

													\section{Preliminaries}
Throughout the text, by a map $f\colon X\to Y$ we always mean a holomorphic map between complex manifolds $X$ and $Y$ of dimensions $n$ and $p$, respectively.

\begin{definition}
Given a map $f\colon X\to Y$, an \emph{unfolding of $f$} over a pointed manifold $(S,s_0)$ is any map $F\colon\cX\to \cY$, endowed with two embeddings $i\colon X\to\cX, j\colon Y\to \cY$ and two submersions $\alpha\colon \cX\to S,\beta\colon \cY\to S$, satisfying:
\begin{enumerate}
\item The following diagram  commutes:

$$\xymatrix{
X\ar[rr]^{f}\ar[d]^{i}&&Y\ar[d]^{j}\\
\cX\ar[rr]^F\ar[dr]_{\alpha}&&\cY\ar[dl]^\beta\\
&S&}$$
\item Let $X_s=\alpha^{-1}(s)$ and $Y_s=\beta^{-1}(s)$, for any $s\in S$. Then  $i$ and $j$ map $X$ and $Y$ isomorphically to $X_{s_0}$ and $Y_{s_0}$. 
\end{enumerate}
  A \emph{local unfolding} (at a subset $A\subseteq X$) is any unfolding $F$ of the restriction of $f$ to some open neighbourhood $U\subseteq X$ of $A$.  
 \end{definition}
For any $s\in S$, we write $f_s=F\vert_{X_s}\colon X_s\to Y_s$.
	\begin{rem}
\begin{enumerate}
\item If $F$ is an unfolding of $f$, then $f$ and $f_{s_0}$ are $\cA$-equivalent. Indeed, every $\cA$-equivalence can be seen as an unfolding where $S$ consists just of one point.
\item If $F$ is an unfolding of $f$ over $(S,s_0)$ and $\sF$ is an unfolding of $F$ over $(T,t_0)$, then $\sF$ is also an unfolding of $f$ over $(S\times T, (s_0,t_0))$.
\item Every unfolding admits the following local form: take local coordinates so that the maps $i,j,\alpha,\beta$ are given by $i(x)=(x,0)$, $j(y)=(y,0)$, $\alpha(s,x)=s$ and $\beta(s,y)=s$. The map $F$ is written in such coordinates as $F(s,x)=(s,f_s(x))$, with $f_0=f$. Hence, the definition of unfolding can be seen as a global coordinate-free version of the definition of unfolding of a multigerm.
\end{enumerate}
	\end{rem}
 
	\begin{definition}\label{defLocalUnfolding}
An \emph{unfolding} of a multigerm $f\colon (\C^n,S)\to (\C^p,0)$ is a multigerm  $$F\colon(\C^r\times\C^n,\{0\}\times S)\to(\C^r\times\C^p,0),$$ of the form $$(s,x)\mapsto(s,f_s(x)),$$ with $f_0(x)=f(x)$. 

Two unfoldings $F$ and $G$ of $f$ are \emph{$\cA$-equivalent as unfoldings} if there exist biholomorphic unfoldings of the identity $\Phi:(\C^r\times\C^n,\{0\}\times S)\to(\C^r\times\C^n,\{0\}\times S)$ and $\Psi:(\C^r\times\C^p,0)\to(\C^r\times\C^p,0)$, such that $\psi\circ F\circ \Phi^{-1}=G$.

An unfolding $F$ of $f$ is called \emph{trivial} if it is $\cA$-equivalent as unfolding to the constant unfolding $\id\times f$. 

A multigerm $f\colon(\C^n,S)\to (\C^p,0)$ is \emph{stable} if every unfolding of $f$ is trivial. 
	\end{definition}

 	\begin{prop}\label{lemExistenceOfStableUnfolding}\label{propExistsStableUnfolding}
Every finite multigerm admits a stable unfolding.
\begin{proof}If $f$ is finite, then it is $\cK$-finite, that is, it has finite singularity type in the sense of Mather. The result follows since any multigerm of finite singularity type admits a stable unfolding \cite[Theorem 2.8]{GibsonWirthmuller76TopologicalStability}.
\end{proof}
	\end{prop}

	\begin{definition}\label{defAStableMap}
We say that a finite map $f\colon X\to Y$ is \emph{stable} (or $\cA$-stable) if, for any $y\in f(X)$, the multigerm of $f$ at $f^{-1}(y)$ is stable.
	\end{definition}

	\begin{definition}\label{defSigma^i(f)}
Given $f\colon X\to Y$, we define $$\Sigma^k(f)=\{x\in X\mid \corank f_x=k\}$$ and $$\widehat\Sigma^k(f)=\bigcup_{i\geq k}\Sigma^i(f).$$
	\end{definition}
	
	\begin{rem}
The set $\widehat\Sigma^k(f)$ may not be the closure of $\Sigma^k(f)$. The map $(x,y)\mapsto (x^2,y^2,xy)$ has no corank 1 points but has a corank 2 point at the origin.
	\end{rem}

 	\begin{lem}\label{lemBoundDimHatSigmaFiniteMap}\label{lemSigmaFiniteMap}
Let $f\colon X\to Y$ be a finite-to-one map and let $n=\dim X$. Then $$\dim \widehat\Sigma^k(f)\leq n-k,$$ for all $1\leq k\leq n$.  As a consequence, $\dim\Sigma^k(f)\leq n-k$.
\begin{proof}
We proceed by induction on $k$: Assume first $\dim\widehat\Sigma^1(f)=n$, then there exists a proper analytic space $Z\subsetneq X$, such that $f$ has constant rank $d<n$ at $X\setminus Z$. For any point $x \in X\setminus Z$, by the constant rank theorem, we can perform some local changes of coordinates in source and target to obtain a map of the form $(x_1,\dots,x_d,0\dots,0)$, which is not finite. This proves the case $k=1$.

Assume $\dim \widehat\Sigma^k(f)\geq n-k+1$. Then, since $\widehat\Sigma^{k}(f)\subseteq\widehat\Sigma^{k-1}(f)$, by induction we have  $\dim \widehat\Sigma^{k}(f)\leq n-k+1$, so the dimension of $\widehat \Sigma^k(f)$ equals $n-k+1$. Let $x$ be a regular point of $\widehat\Sigma^{k}(f)$ where the dimension of $\widehat\Sigma^{k}(f)$ is $n-k+1$. Then,  there exists an open neighbohood $U\subseteq X$ of $x$, such that the restriction $g=f\vert_{\widehat\Sigma^k(f)\cap U}$ is a holomorphic map defined on a manifold of dimension $n-k+1$. Being a restriction of $f$, the map $g$ is finite-to-one. Since the restriction is done precisely at $\widehat\Sigma^k(f)$, it is obvious that $g$ has rank $\leq n-k$ at all source points. This means $\dim\widehat \Sigma^1(g)=n-k+1$, contradicting the induction hypothesis.
\end{proof}
	\end{lem}

	\begin{ex}
The previous bound is sharp, since it is exact for the map $\C^n\to\C^p$ given by
$$(x_1,\dots,x_n)\mapsto (x_1^2,\dots,x_n^2,0,\dots,0).$$
	\end{ex}

The following fact is well known (see for instance \cite{GibsonSingularPointsSmoothMappings}).

	\begin{prop}\label{propDimSigmaStableMap}
For any stable map $f\colon X\to Y$, the space $\Sigma^k(f)$ is empty or a manifold of dimension $n-k(p-n+k)$. 
 	\end{prop}
\subsection*{Notation on products of copies of $X$:}We write elements in $X^k$ as tuples $w=(x^{(1)},\dots,x^{(k)})$ of points $x^{(l)}\in X$, each one with local coordinates $x^{(l)}_1,\dots,x^{(l)}_n$. For $k=2,3$, we use $x=x^{(1)},x'=x^{(2)},x''=x^{(3)}$ and denote the coordinates by $x_i=x^{(1)}_i,x'_i=x^{(2)}_i$ and $x''_i=x^{(3)}_i$. 
The \emph{small diagonal} of $X^k$ is the subset $$\Delta(X,k)=\{(x^{(1)},\dots,x^{(k)})\in X^k\mid x^{(1)}=\dots=x^{(k)}\}.$$ The \emph{big diagonal} of $X^k$ is the subset $$D(X,k)=\{(x^{(1)},\dots,x^{(k)})\in X^k\mid x^{(i)}=x^{(j)}\text{ for some $i\ne j$}\}.$$
When working locally, we write $\Delta(n,k)$ and $D(n,k)$ for the germs of $\Delta(\C^n,k)$ and $D(\C^n,k)$ at 0. The ideals in $\cO_{kn}$ defining these two space germs are 
\begin{align*}
I_{\Delta(n,k)}&=\sum_{l=2}^k\langle x^{(1)}_i-x^{(l)}_i\mid 1\leq i\leq n\rangle,\\
I_{D(n,k)}&=\bigcap_{1\leq l<m\leq k}\langle x^{(l)}_i-x^{(m)}_i\mid 1\leq i\leq n \rangle.
\end{align*}

	\begin{definition}\label{defNormalCrossings}
We say that a map $f\colon X\to Y$ has \emph{normal crossings} if, for any $k\geq 2$, the restriction of $f^k$ to $X^k\setminus D(X,k)$ is transverse to $\Delta(Y,k)$.
	\end{definition}
	
	An easy transversality argument shows the following well known result:
		\begin{prop}\label{propStableImpliesNormalCrossings}
Any stable map has normal crossings.
	\end{prop}

													\section{Multiple Points}\label{M1M2}

	\begin{definition}\label{defStrictMultPoint}
	Given a map $f\colon X\to Y$, we say that $(x^{(1)}, \dots, x^{(k)})\in X^k$ is a \emph{strict $k$-multiple point of $f$} if $f(x^{(i)})=f(x^{(j)})$ and $x^{(i)}\neq x^{(j)}$, for all $i\neq j$. We denote by $D^k_S(f)$  the analytic closure of the set of strict $k$-multiple points of $f$, that is, 
	$$
	D^k_S(f)=\overline{(f^k)^{-1}(\Delta(Y,k))\setminus D(X,k)}.
	$$
We regard $D^k_S(f)$ as a complex space with the reduced structure. If we write $\sI_{\Delta(X,k)}$ and $\sI_{D(X,k)}$  for the ideal sheaves in $X^k$ of the small and big diagonal respectively, then (see the beginning of Section \ref{secAnotherMultiplePointStructure}) the ideal sheaf of $D^k(f)$ is 
$$
\sqrt{\sP(f,k): {\sI_{D(X,k)}}^\infty},
$$
where,  $$\sP(f,k)=(f^k)^*\sI_{\Delta(Y,k)}.$$
	\end{definition}

We have local versions of the above definitions as well: For any finite multigerm $f\colon (X,A)\to (Y,y)$, we take a  representative $\hat f$ of $f$, defined at a small enough neighbourhood $U$ of $A$. Then, we define $D^k_S(f)$ as the multigerm at $A^k$ of $D^k_S(\hat f)$. 
As in the global case, we have $D^k_S(f)=\overline{(f^k)^{-1}(\Delta(p,k))\setminus D(n,k)},$
 whose defining ideal in $\cO_{kn}$ is $\sqrt{P(f,k): {I_{D(n,k)}}^\infty},$ with $P(f,k)=(f^k)^*I_{\Delta(p,k)}$.


	\begin{ex}\label{exDeformationOfCusp}
Take the family of curves $f_t\colon \C\to \C^2$, given by 
$$x\mapsto (x^2,x^3+tx).$$
The map $f_0$ is just a usual cusp. It is injective, and thus $D^k_S(f_0)$ is empty. For $t\neq 0$, a straightforward computation shows that $D^k_S(f_t)$ is the zeroset of $\langle x^2+t,x+x'\rangle$.  
	\begin{figure}[ht]
\begin{center}
\includegraphics[scale=0.85]{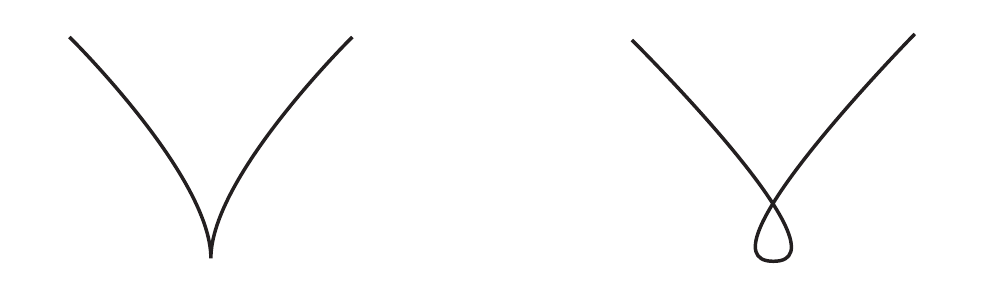}
\end{center}
\caption{Images of $f_0$ and $f_t,t\neq 0$, respectively}
\end{figure}

Let $I=\langle x^2+t,x+x'\rangle$ and denote by $I_{t_0}$ the ideal in $\cO_2$ (variables $x$ and $x'$) obtained by the substitution $t=t_0$ in $I$.  Then:
\begin{itemize}
\item $I_0$ does not define $D^k_S(f_0)$.
\item $\cO_2/I_0$ is not reduced.
\end{itemize}
	\end{ex}

As the previous example shows, the space $D^k_S(f)$ does not behave well under deformations. To avoid this problem, we need to adopt a different definition of the multiple point spaces, including some non-strict multiple points and allowing the space to be non-reduced. 
Before getting into the details, we show the following:

 	\begin{lem}\label{propiedades}
For any multigerm $f\colon (\C^n,A)\to(\C^p,0)$.
\begin{enumerate}
\item If $F=(s,f_s)$ and $F'=(s,f'_{s})$ are $\cA$-equivalent as unfoldings of $f$, then $$D^k_S(F)\cap \{s=0\}=D^k_S(F')\cap \{s=0\}.$$
\item If $f$ is stable, then for any  unfolding $F=(s,f_s)$ of $f$, 
$$D^k_S(F)\cap\{s=0\}=D^k_S(f).$$
 \item If $F=(s,f_s)$ and $F'=(t,f'_{t})$ are stable unfoldings of $f$,
then 
$$D^k_S(F)\cap \{s=0\}=D^k_S(F')\cap \{t=0\}.$$
\item Let $F=(s,f_s)$ and $F'(t,f'_t)$ be stable unfoldings of two germs $f,f'$, respectively. If $f'=\psi\circ f\circ \phi^{-1}$, where $\phi,\psi$ are biholomorphisms, then 
$$
\phi^k(D^k_S(F)\cap\{s=0\})=D^k_S(F')\cap\{t=0\}.
$$
\end{enumerate}
\begin{proof} 

\medskip
1) By hypothesis, we have $F'=\Psi \circ F\circ \Phi^{-1}$, where $\Phi,\Psi$ are unfoldings of the identity in $\C^n,\C^p$ respectively. On one hand, $\Phi(D^k_S(F))=D_S^k(F')$. On the other hand, writting $\Phi=(s,\phi_s)$  with $\phi_0=\id$,
	$$
	D^k_S(F)\cap\{s=0\}=\Phi(D_S^k(F))\cap\{s=0\}=D_S^k(F')\cap\{s=0\}.
	$$
	
2)  Since $f$ is a stable map and $F$ is an unfolding of $f$, then $F$ is $\cA$-equivalent to the constant unfolding $(\id\times f)$. From (1), it follows 
$$
D^k_S(F)\cap\{s=0\}=D^k_S(\id\times f)\cap\{s=0\}=D^k_S(f).
$$

3) Let $\cF$ be the germ given by $\cF(s,t,x)=(s,t,f_s(x)+f'_t(x)-f(x))$. This is an unfolding of both $F$ and $F'$. Since $F$ and $F'$ are stable, (2) implies 
$$D^k_S(F)\cap\{s=0\}=D^k_S(\cF)\cap\{s=t=0\}=D^k_S(F')\cap\{t=0\}.
$$

4) Let $G$ be the unfolding of $f'$ given by $G(s,x)=(s,\psi\circ f_s\circ\phi^{-1})$. We have that $\Psi\circ F\circ \Phi^{-1}=G$, where $\Phi=\id\times\phi$ and $\Psi=\id\times\psi$, hence $G$ is also stable. From (3) we obtain
$$
D^k_S(G)\cap\{s=0\}=D^k_S(F')\cap\{t=0\}.
$$
On the other hand, $\Phi^k(D^k_S(F))=D^k_S(G)$ and thus,
$$
\phi^k(D^k_S(F)\cap\{s=0\})=\Phi^k(D^k_S(F))\cap\{s=0\}=D^k_S(G)\cap\{s=0\}.
$$	
\end{proof}

	\end{lem}

 Now we show that there is a unique way to define $k$th multiple point subspaces of $X^k$ satisfying some natural requirements.  To be precise about these scheme-theoretic requirements, we need some definitions: A \emph{multiple point space structure} is a rule, denoted by $R^k$, which associates to any finite map $f\colon X\to Y$ a closed complex subspace $R^k(f)$ of $X^k$. Since we ask $R^k(f)$ to be a closed subspace of $X^k$, the structure sheaf $\cO_{R^{k}(f)}$ is defined by some coherent ideal sheaf $\sI_R^k(f)$ in $\cO_{X^k}$.

 	\begin{definition}\label{defConditionsM1M2}
We define the following two conditions for any multiple point space structure:
\begin{itemize}
\item[M1.] If $f$ is a stable finite map, then $R^k(f)=D^k_S(f)$.
\item[M2.] For any local unfolding $F$ of $f$ at an open neighbourhood  $U\subseteq X$, the map $i^k$ sends $R^k(f)\cap U^k$ isomorphically to $R^k(F)\cap (U_{s_0})^k$.
\end{itemize}
	\end{definition}

Condition M1 may be thought as $R^k$ being the simplest choice for the multiple point scheme. Condition M2 means, first, that $R^k$ behaves well under deformations and, second, that $R^k(f)$ is determined by the multilocal behaviour of $f$ at every collection of points. Now we show that these two conditions determine $R^k$ uniquely.

	\begin{defprop}\label{propUniqueMultiplePointStructure}
There exists a unique multiple point structure $D^k$ satisfying M1 and M2.  For any finite map $f\colon X\to Y$, we call $D^k(f)$ the \emph{$k$th multiple point space} of $f$. For any point $w\in X^k$, the space $D^k(f)$ is given locally around $w$ by
$$
D^k(f)=(i^k)^{-1}(D^k_S(F)\cap (X_{s_0})^k),
$$ 
where $F$ is any local stable unfolding of $f$ at $w$. 
\begin{proof}
First we explain how $D^k(f)$ is defined. We set $$D^k(f)\cap(X^k\setminus(f^k)^{-1}(\Delta(Y,k)))=\emptyset.$$ Let $x^{(1)},\dots,x^{(k)}\in (f^k)^{-1}(\Delta(Y,k))$ and take $f\colon (X,A)\to (Y,y)$, where $A=\{x^{(1)},\dots,x^{(k)}\}$ and $f(x^{(i)})=y$. We define $D^k(f)$ in a neighbourhood of $(x^{(1)},\dots,x^{(k)})$ as
$$
D^k(f)=\phi^k(D_S^k(F)\cap\{s=0\}),
$$
where $F$ is any stable unfolding of $f$. By Lemma \ref{propiedades}, $D^k(f)$ is well defined and does not depend on the choice of $F$. Therefore, these spaces can be glued together to get a  complex space defined globally. We will write it as $D^k(f)$ and its defining ideal sheaf as $\sI^k(f)$. It follows from the definition that $D^k(f)$ is given by 
$$
D^k(f)=(i^k)^{-1}(D^k_S(F)\cap (X_{s_0})^k),
$$ 
where $F$ is any local stable unfolding of $f$.

Now we show that the construction satisfies M1 and M2. If $f$ is stable, then we can take $F=f$ and hence $D^k(f)=D^k_S(f)$, so condition M1 is satisfied. Condition M2 is  obvious as well. Taking local coordinates, it suffices to prove the claim for a multigerm  $f\colon (\C^n,A)\to(\C^p,0)$. Given any unfolding $F=(t,f_t)$ of $f$, we take $\cF(s,t,f_{s,t})$ a stable unfolding of $F$. Then,
\begin{align*}
D^k(f)&=D^k_S(\cF)\cap\{s=t=0\}=(D^k_S(\cF)\cap\{s=0\})\cap\{t=0\}\\
&=D^k(F)\cap\{t=0\}.
\end{align*}

Finally, we show the unicity. Let $\hat D^k$ be another $k$th multiple point structure satisfying M1 and M2. For any map $f$ and for any local stable unfolding $F$ of $f$, we have locally the following equalities:
$$
D^k(f)=(i^k)^{-1}(D^k_S(F)\cap (X_{s_0})^k)=(i^k)^{-1}(\hat D^k(F)\cap (X_{s_0})^k)=\hat D^k(f).
$$
\end{proof}

	\end{defprop}

	\begin{ex}
We compute the double point space of the cusp $f\colon (\C,0)\to (\C^2,0)$ given by $$x\mapsto (x^2,x^3).$$Take the family of curves of Example \ref{exDeformationOfCusp} as an unfolding, that is, take the map germ $F\colon(\C^2,0)\to(\C^3,0)$, given by 
$$(t,x)\mapsto(t,x^2,x^3+tx).$$

Let $\phi$ and $\psi$ be the local changes of coordinates $(t,x)\mapsto(t-x^2,x)$ and $(X,Y,Z)\mapsto (X-Y^2,Y,Z)$. The map $\phi\circ F\circ\psi$ is a cross-cap, and therefore $F$ is stable. The strict double point space $D^2_S(F)$ is defined by the ideal $$\sqrt{P(F,2): {I_{D(n,2)}}^\infty}=\langle t-t',x+x',t+x^2\rangle.$$
The substitution $t=0$ yields the double point space
 $$D^2(f)=V(x+x',x^2).$$
\begin{figure}[ht]
\begin{center}
\includegraphics[scale=0.8]{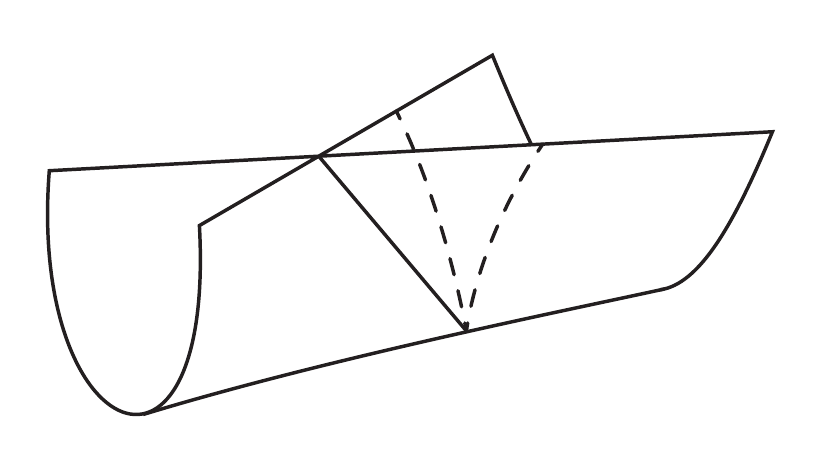}
\caption{Image of the cross-cap, as an unfolding of the cusp.}
\end{center}
\end{figure}
	\end{ex}

	\begin{ex}\label{exAMuHasVeryBigStableUnfolding}
Let $A_\mu\colon (\C,0)\to (\C,0)$ be the germ given by $x\mapsto x^{\mu+1}.$ To compute its multiple point spaces, we need to take a stable unfolding $F$ of $A_\mu$. One can check (see \cite{GibsonSingularPointsSmoothMappings} for details) that the map germ $F\colon(\C^{\mu-1}\times\C,0)\to (\C^{\mu-1}\times\C,0)$, given by $$(u_1,\dots, u_{\mu-1},x)\mapsto (u_1,\dots, u_{\mu-1},x^{\mu+1}+u_{\mu-1}x^{\mu-1}+\dots+u_1x),$$
is a minimal (in the sense of the number of parameters $u_i$) stable unfolding of $A_\mu$. 
	\end{ex}
	
As the previous example shows, the construction of the multiple point space forces us to study maps between manifolds of arbitrarily big dimension. For double points in any corank and for $k$-multiple points of corank 1 map germs, this problem is adressed by Mond's ideals (Section \ref{secMultiplePointsCorank1} and Theorem \ref{thmI2DefinesD2}).
	
	\begin{prop}\label{lemEliminationParametersRank}
Let $f\colon(\C^n,0)\to(\C^p,0)$ be a rank $r$ finite map germ. 

\begin{enumerate}
\item If $f$ is of the form $(s,x)\mapsto(s,f_s(x)), s\in \C^r, x\in \C^{n-r}$, then the projection $P\colon \C^{nk}\to\C^{r}\times\C^{k(n-r)}$ which forgets the variables $s^{(2)},\dots,s^{(k)}$ induces an isomorphism
$$D_S^k(f)\cong\overline{\{(s,w)\in \C^r\times \C^{k(n-r)}\mid w\text{ is a strict multiple point of }f_s\}}.$$
\item $D^k(f)$ embeeds into $\C^r\times\C^{k(n-r)}$.
\end{enumerate}

\begin{proof}
1) Set $Z=\{(s,w)\in \C^r\times \C^{kn}\mid w\text{ is a strict multiple point of }f_s\}$ and observe that 1) is a set theoretical question, since both $D^k_S(f)$ and $\overline{Z}$ are reduced. It is obvious that $P$ restricts to a bijection $$\{\text{strict $k$-multiple points of } f\} \to Z.$$ Therefore, $P(D^k_S(f))\subseteq \overline Z$. Let $\gamma\colon D\to \overline{Z}$ be a curve defined in a neighboohood of the origin of $\C$, satisfying $\gamma(D\setminus\{0\})\subseteq Z$. Let $\gamma_ i$ be the coordinate functions of $\gamma$ and  let $$\sigma=(\gamma_1,\dots,\gamma_r)\text{ and }\omega=(\gamma_{r+1},\dots,\gamma_{r+k(n-r)}).$$ Let $\gamma'\colon D\to X^k$ be defined by $$\gamma'(t)=(\sigma(t),\dots,\sigma(t), \omega(t)),$$
with $\sigma$ repeated $k$ times. It is obvious that $\gamma'(t)$ is a strict $k$-multiple point of $f$, for all $t\in D\setminus\{0\}$. Therefore, we have $\gamma'(D)\subseteq D^k_S(f)$. Since $P\circ \gamma'=\gamma$, we obtain $\overline Z\subseteq P(D^k_S(f))$, as desired.

To show 2) take $F\colon (\C^l\times\C^n,0)\to (\C^l\times\C^p)$ a stable unfolding of $f$ of the form $F(t,s,x)=(t,s,f_{t,s}(x))$. By 1), we can eliminate the variables $t^{(2)},\dots,t^{(k)}$ and $s^{(2)},\dots,^{(k)}$ to obtain an isomorphic image of $D^k_S(F)$ into $\C^{l+r}\times\C^{k(n-r)}$. Now, since $D^k(f)=D^k_S(F)\cap\{t^{(i)}=0\mid 1\leq i\leq k\}$, the claim follows immediately.
\end{proof}

	\end{prop}

	\begin{lem}\label{lemDimensionStrictStable}
 If $f\colon X\to Y$ is a stable map, then the set of strict $k$-multiple points of $f$ is empty or a manifold of dimension $kn-(k-1)p$. In particular, $D^k(f)$ is empty or reduced and equidimensional of dimension $kn-(k-1)p$.
 \begin{proof}The statement follows directly from Proposition \ref{propStableImpliesNormalCrossings} , taking into account that $\Delta(Y,k)$ is a manifold of codimension $(k-1)p$ in $Y^k$.\end{proof}
	\end{lem}

	\begin{prop}\label{lemDimensionBoundMultiplePoints}
For any finite map $f\colon X\to Y$, any integer $k\geq 2$ and any point $w\in D^k(f)$, all non embedded irreducible components of the germ $(D^k(f),w)$ have dimension $\geq kn-(k-1)p$.
In particular, the $k$th multiple point space $D^k(f)$ is empty or has dimension  $\geq kn-(k-1)p$. 
\begin{proof}
Let $F(s,x)=(s,f_s(x)), s\in \C^r$ be a local stable unfolding of $f$, so that $D^k(f)=D^k(F)\cap \{s=0\}$. By Lemma \ref{lemDimensionStrictStable}, $D^k(F)$ is empty or has  dimension $k(r+n)-(k-1)(r+p)=kn-(k-1)p+r$. Let $w \in D^k(f)$, then  $(0,w)\in D^k(F)$. By Lemma \ref{lemDimensionStrictStable}, all the irreducible components of $D^k(F)$ have dimension exactly $kn-(k-1)p+r$. The statement follows since $D^k(f)$ is obtained by intersecting $D^k(F)$ with $\{s=0\}$, which is a manifold of dimension $r$. 
\end{proof}
	\end{prop}
	
	\begin{definition}\label{defDimensionallyCorrect}
For any finite map $f\colon X\to Y$  the $k$th multiple point space  $D^k(f)$ is \emph{dimensionally correct} if it is empty or has dimension $kn-(k-1)p$.
	\end{definition}
	
\begin{open}
 Is $D^k(f)$ Cohen Macaulay if it is dimensionally correct? For corank 1 map germs the answer is positive, as we will see in Section \ref{secMultiplePointsCorank1}. For arbitrary corank, Theorem \ref{DimCorrectImpliesD2(f)CM} gives a positive answer, but only for double points.
\end{open}

								\section{Multiple points of corank 1 monogerms}\label{secMultiplePointsCorank1}

								\label{LevCorrango1}
In the previous section we have shown that there is a unique multiple point structure satisfying some reasonable conditions. We have also shown how to compute these structure by taking a stable unfolding and slicing the closure of its strict $k$-multiple points. However, the computation of the multiple point spaces can be hard in practice, as Example \ref{exAMuHasVeryBigStableUnfolding} shows.

In \cite{MondSomeRemarks} Mond introduces some ideals $I^k(f)$ for corank 1 map germs which define the multiple point spaces. These ideals are obtained directly from the original map, with no need to take any unfolding. Moreover, in \cite{MararMondCorank1} Marar and Mond  show that stability and finite determinacy of corank 1 map germs $f\colon (\C^n,0)\to(\C^p,0),\ n<p,$ can be characterized by the
geometry of the multiple point spaces.

 Let $f\colon (\C^n,0)\to(\C^p,0)$ be a  corank $1$ map germ with $n\leq p.$ Up to $\mathcal
A$-equivalence, $f$ can be written in the form $$(x,y)\mapsto(x,f_n(x,y),\dots,f_p(x,y)),$$ with $x\in\C^{n-1}$
and $y\in\C.$ We can think of $f(x,y)$ as a $(n-1)$-parameter family of
functions of one variable $f_x(y)=(f_n(x,y),\dots,f_p(x,y))$. Embedding $D^2(f)$ in $\C^{n-1}\times \C^2$ (see Proposition \ref{lemEliminationParametersRank}), a point $(x,y,y')$ is a double point if and only if the coefficients of the Newton interpolating polynomial of degree $1$ for
the points $(y,f_{j,x}(y)),(y',f_{j,x}(y'))$ are equal to $0,$ for all $n\leq j\leq p$.
These coefficients, the generators of \emph{Mond's double point ideal} $I^2(f)$, are the divided differences $$f_{j,x}[y,y']=\frac{f_j(x,y)-f_j(x,y')}{y-y'}.$$ Similarly, the
triple points are $(x,y,y',y'')$ $\in\C^{n-1}\times \C^3$ such that every coefficient of the Newton
interpolating polynomial of degree $2$ for the points $(y,f_{j,x}(y)),(y',f_{j,x}(y')),(y'',f_{j,x}(y''))$ are
equal to $0$ for every $n\leq j \leq p$. These coefficients are the divided differences $f_{j,x}[y,y']$ and the iterated divided differences
$$f_{j,x}[y,y',y'']=\frac{f_{j,x}[y,y']-f_{j,x}[y,y'']}{y'-y''}.$$ Hence, \emph{Mond's triple point ideal} is
$${I}^3(f)=\langle f_{j,x}[y,y'],f_{j,x}[y,y',y''] \mid n\leq j\leq p\rangle.$$ Higher multiple ideals $I^k(f)$ are defined
analogously.

	\begin{rem}\label{remLagrangeInterpolation}
As observed by Marar and Mond \cite{MararMondCorank1}, the coefficients of the Lagrange interpolation polynomial provide another set of generators for the ideal
defining the $k$-multiple points, with the advantage that they are invariant under the action of the symmetric
group $S_k.$
	\end{rem}

	\begin{prop}\cite[2.16]{MararMondCorank1}\label{propGaffneyMondAgreeCorank1}
For any finite corank 1 map germ $f\colon(\C^n,0)\to(\C^p,0)$, the ideal $I^k(f)$ defines $D^k(f)$. In particular, $D^k(f)$ is a complete intersection, if it is dimensionally correct.
	\end{prop}

\begin{thm}\label{Marar-MondEstableCor1}\cite[Thm. 2.14]{MararMondCorank1}
Let $f\colon(\C^n,0)\to (\C^p,0)$ be a finite corank 1 map, with $n<p$. Then: 
\begin{enumerate}
\item $f$ is stable if and only if ${D}^k(f)$ is empty or smooth of dimension $kn-(k-1)p$, for every $k\geq 2$.
\item $f$ is finitely determined if and only if
 ${D}^k(f)$ is empty or an ICIS of dimension $kn-(k-1)p$, for all $k$ satisfying
$kn-(k-1)p\ge0$, and ${D}^k(f)\subseteq\{0\}$ for $kn-(k-1)p<0$.
   \end{enumerate}
\end{thm}

\begin{rem}
As we will show in Proposition \ref{D2NotRegularInCorank2}, the double point locus of a stable map may contain singularities, as long as the map has some corank $\geq 2$ points. Thus, any stability or finite determinacy criterion in terms of $D^k(f)$ for corank $\geq 2$ map germs must necessarily be more complicated than the one above. 
\end{rem}

								\section{The double point ideal sheaf}\label{secDOublePointIdealSheaf}

In  \cite{MondSomeRemarks}, Mond gives a explicit set of generators for a double point ideal of any map germ $f:(\C^n,0)\to(\C^p,0)$.  The aim of this section is to show (1) how these local constructions glue to global schemes for any map, and (2) that the double point structure obtained satisfies conditions M1 and M2, and thus it coincides with the double point structure of Section \ref{M1M2}.

	\begin{defprop}\label{defDoublePointSpaceOfGerm}
For any map germ $f\colon(\C^n,0)\to(\C^p,0)$, the germs $f_j(x)-f_j(x'), 1\leq j\leq p$ vanish on the diagonal $\Delta(n,2)$. Therefore, they are contained in the ideal generated by $x_i-x'_i, 1\leq i\leq n$. In other words, for all $j\leq p$, there exist some function germs $\alpha_{ij}\in \cO_{2n}$, satisfying 
$$f_j(x)-f_j(x')=\sum_{i=1}^n\alpha_{ji}(x,x')(x_i-x_i').
$$ 
This can be expressed as the matrix equality
$$f(x)-f(x')=\alpha (x-x'),$$
 where  $\alpha$ represents the $p\times n$ matrix $(\alpha_{ji})$ and $x-x'$ and $f(x)-f(x')$ are taken as column vectors of sizes $n$ and $p$ respectively.  
\emph{Mond's double point ideal} $I^2(f)$ is the sum 
$$I^2(f)=(f\times f)^*(I_{\Delta(p,2)})+\langle \text {$n\times n$ minors of }\alpha\rangle.
$$ 
The matrix $\alpha$ may not be unique, but $I^2(f)$ does not depend on the choice of $\alpha$ \cite[Prop 3.1]{MondSomeRemarks}. 
	\end{defprop}
	
	\begin{ex}\label{exDoubleFoldTrifold}
Let $f\colon(\C^2,0)\to(\C^3,0)$ be given by
$$(x,y)\mapsto(x^2,y^2,x^3+y^3+xy).$$
\begin{figure}[ht]
\begin{center}
\includegraphics[scale=0.7]{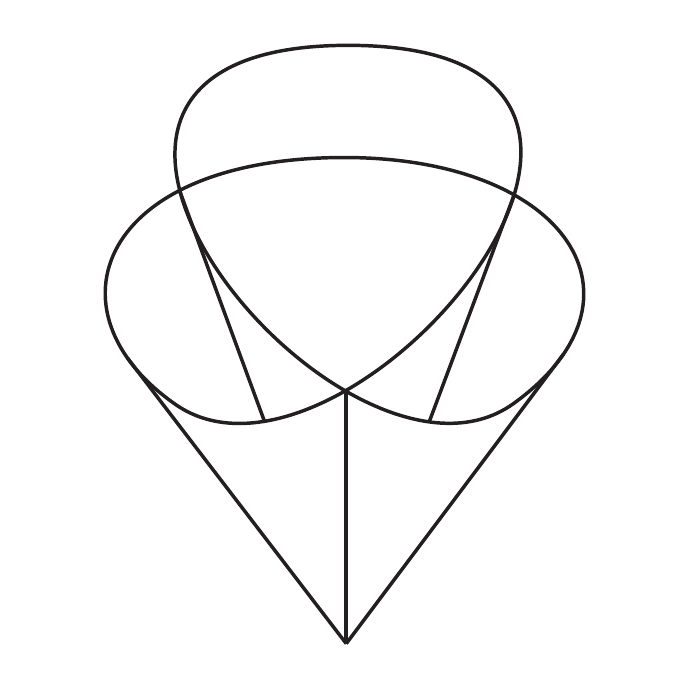}
\caption{Image of the map germ in Example \ref{exDoubleFoldTrifold}.}
\end{center}
\label{figDoubleFoldTrifold}
\end{figure}

Computing the double point space of $f$ as explained in Proposition-Definition \ref{propUniqueMultiplePointStructure} can be quite difficult. However, Mond's ideal $I^2(f)$ can be computed easily: A  solution for the equation $f(x)-f(x')=\alpha (x-x')$
 is 
$$\left(\begin{array}{c}x^2-x'^2 \\y^2-y'^2 \\x^3+y^3+xy-x'^3-y'^3-x'y'\end{array}\right)=$$
$$=\left(\begin{array}{cc}x+x' & 0 \\0 & y+y' \\x^2+xx'+x'^2+y & y^2+yy'+y'^2+x'\end{array}\right)\left(\begin{array}{c}x-x' \\y-y'\end{array}\right).$$
Now $I^2(f)$ is generated in $\cO_{2n}$ by
\begin{align*}
g_1&=(x+x')(x-x'), &g_4&=(y + y') (2 x^2 + 2 x x' + 2 x'^2 + y + y'),\\
g_2&=(y+y')(y-y'),  &g_5&=(x + x') (x + x' + 2 y^2 + 2 y y' + 2 y'^2),\\
g_3&=(x+x')(y+y'),  &g_6&=(x - x') (2 x^2 + 2 x x' + 2 x'^2 + y + y')+\\
& &&+(y - y') (x + x' + 2 y^2 + 2 y y' + 2 y'^2).
\end{align*}

	\end{ex}

%
%
%

	\begin{lem}\label{lemAlphaAgreesDifferentialOnDiagonal}
Let $f\colon (\C^n,0)\to (\C^p,0)$ be a map germ and let $\alpha $ be as above, then $\alpha(x,x)=df_x.$ 
\begin{proof}Let $e_i$ be the $i$-th vector of the canonical basis of $\C^n$. Then
	$$\alpha_{ji}(x,x)=\lim_{\lambda\to0}\alpha_{ji}(x,x+\lambda e_i)=\lim_{\lambda\to0}\frac{f_j(x)-f_j(x+\lambda e_i)}	
	{\lambda}=\pd{f_j}{x_i}.$$

\end{proof}

	\end{lem}

Next lemma, which can be shown by elementary techniques, shows that $I^2(f)$ behaves well under under $\cA$-equivalence.

	\begin{lem}\label{lemI2AEquivalenceIsomorphism}
Let $f$ and $g$ be $\cA$-equivalent map germs with $f=\psi\circ g\circ\varphi$. Then, $(\varphi\times\varphi)^*(I^2(f))=I^2(g)$.
\end{lem}

%

Given a map $f\colon X\to Y$ and a point $x\in X$, we take local coordinates so that the germ $f_x$ of $f$ at $x$ is $f_x=\psi\circ f'\circ \varphi^{-1}$, for some biholomorphisms $\varphi,\psi$ and some map germ $f'\colon (\C^n,0)\to(\C^p,0)$. We define the ideal $I^2(f_x)$ in $\cO_{X\times X,(x,x)}$ as
$$
I^2(f_x)=(\varphi\times\varphi)^*(I^2(f')).
$$
Lemma \ref{lemI2AEquivalenceIsomorphism} ensures that this definition does not depend on the choice of $\varphi$ and $\psi$. The following lemma allows us to extend the local definition of $I^2(f)$ to a global ideal sheaf of double points. 
 
 	\begin{lem}\label{lemStalksOfMondStructure}
Given a map $f\colon X\to Y$ and a point $x\in X$, denote by $f_x$ the germ of $f$ at $x$. Let $\sI^2(f)$ be the ideal sheaf given by some representatives of the generators $I^2(f_{x})$ on a sufficiently small open neighbourhood $U$ of $(x,x)$ in $X\times X$. The following hold:
\begin{enumerate}
\item $\sI^2(f)_{(x',x')}=I^2(f_{x'})$, for any $(x',x')\in \Delta(U,2)$.
\item $\sI^2(f)_{(x',x'')}=\sP(f,2)_{(x',x'')}$, for any $(x',x'')\in U\setminus \Delta(U,2)$.
\end{enumerate}
\begin{proof}
By Lemma \ref{lemI2AEquivalenceIsomorphism} we can assume $X=\C^n$ and $Y=\C^p$, taking local coordinates.
To show (1), we just need to shrink $U$ so that we have representatives of the germs $f_j(x)-f_j(x')$ and of the entries of the matrix $\alpha$ defined on all $U$. Therefore, the germs of these representatives at $(x',x')$ produce the corresponding germs and the corresponding matrix around $(x',x')$. 

To show (2) we need to show that the ideal generated by the germs at $(x',x'')$ of the $n\times n$ minors of $\alpha$ is contained in $\sP(f,2)_{(x',x'')}$. Let $A$ be the submatrix of $\alpha$ obtained by picking the rows $j_1,\dots, j_n$ of $\alpha$. Let $b$ be the vector with entries $f_{j_1}(x)-f_{j_1}(x'),\dots, f_{j_n}(x)-f_{j_n}(x')$. Since $(x',x'')$ is not a diagonal point, there exists $i\leq n$ such that $x'_i\neq x''_i$. Let $A'$ be the matrix obtained by substitution of the $i$-th column of $A$ by $b$. By Cramer's Rule we obtain $\vert A\vert=\vert A'\vert/(x'_i-x''_i)\in \sP(f,2)_{(x',x'')}$ 
\end{proof}
	\end{lem}
	
	The second item of the previous lemma is equivalent to \cite[Lemma 2.1.17]{AltintasThesis}

\begin{definition} The \emph{sheaf of double points} $\sI^2(f)$ of a map $f\colon X\to Y$ is defined as the glueing of the following local structures: Off the diagonal, $\sI^2(f)$ is just the restriction of the sheaf $\sP(f,k)$ to $X^2\setminus D(X,2)$. If $(x,x)$ is a diagonal point, then at a neighbourhood of $(x,x)$ the sheaf is locally given by the double point ideal $I^2(f)$ of the germ of $f$ centered at $x$.
\end{definition}
To glue these local structures, we need to check that, if we compute the structure locally around some point, then the stalk of this local structure at any other close enough point agrees with the structure computed at this other point. This is precisely Lemma \ref{lemStalksOfMondStructure}.

	\begin{lem}\label{lemZerosetI2}
Set theoretically,  $V(\sI^2(f))$ is the union of the strict double points of $f$ and the pairs $(x,x)$ such that  $f$ is singular at $x$.
\begin{proof} Let $(x,x'), x\neq x'$ a non diagonal point in $X^2$. Locally, $\sI^2(f)$ equals $\sP(f,2)$, which vanishes if and only if $(x,x')$ is a strict double point of $f$. Let $(x,x)$ be a diagonal point in $X\times X$, then $\sP(f,2)$ vanishes trivially at $(x,x)$. Moreover, if $\alpha$ is the matrix in the definition of $I^2(f_x)$, by Lemma \ref{lemAlphaAgreesDifferentialOnDiagonal} $\alpha(x,x)$ equals the differential matrix of $f$ at $x$. Therefore, the $n \times n$ minors of $\alpha$ vanish at $x$ if and only if $f$ is singular at $x$.
\end{proof}
	\end{lem}

The following lemma, which appears in \cite{Bivia-Nuno}, can be obtained easily from results about Cohen-Macaulay modules (see \cite{MatsumuraRingTheory} for details).

	\begin{lem}\label{lemCMPullbackCodim}
Let $\phi\colon (\C^m,0)\to (\C^r,0)$ be any map germ. Let $I$ be an ideal in $\cO_r$ and let $J=\phi^*(I)$. If $\cO_r/I$ is Cohen-Macaulay and $\codim V(I)=\codim V(J)$, then $\cO_m/J$ is Cohen-Macaulay. 
	\end{lem}
With a trivial modification of the proof of \cite[Prop. 2.1.11]{AltintasThesis}, we obtain
	\begin{lem}\label{lemConcciniStrickland}
Let $f\colon (\C^n,0)\to (\C^p,0)$ with $n\leq p$. Then
\begin{enumerate}
\item If $\cO_{2n}/I^2(f)\neq0$ then $\dim\cO_{2n}/I^2(f)\geq 2n-p$.
\item If $\dim\cO_{2n}/I^2(f)=2n-p$, then $\cO_{2n}/I^2(f)$ is Cohen Macaulay.
\end{enumerate}
\begin{proof}
We identify the space of $n\times p$ matrices $A=(a_{ji})$ and $n\times 1$ vectors $(d_1,\dots,d_n)^T$ with $\C^{np}\times\C^n$. Let $I$ be the ideal in $\cO_{np+n}$ generated by the entries of $Ad$ and the $n\times n$ minors of $A$. Set $D=V(I)\subseteq \C^{np}\times\C^n$. It turns out that $D$ is  a Buchsbaum-Eisenbud variety of complexes (more precisely $D=W(n-1,1)$, with $n_0=p, n_1=n$ and $n_2=1$, in the notation of \cite{ConciniStricklandOntheVarietyOfComplexes}). By \cite[Thm 2.7, Lemma 2.3]{ConciniStricklandOntheVarietyOfComplexes} $D$ is a Cohen Macaulay subspace of  $\C^{np}\times\C^n$ of codimension $p$.

Now given a matrix $\alpha$ satisfying $f(x)-f(x')=\alpha(x-x')$, we take the map $\phi\colon \C^n\times \C^n\to \C^{np}\times \C^n$ given by
$$(x,x')\mapsto (\alpha(x,x'),(x_1-x'_1,\dots, x_n-x'_n)).$$
We obtain $I^2(f)=\phi^*(I)$, and the result follows directly from Lemma \ref{lemCMPullbackCodim}
\end{proof}
	\end{lem}

We can patch the local results above to obtain the following.

	\begin{thm}\label{thmD2CohenMacaulay}
For any map $f\colon X\to Y$, if $V(\sI^2(f))$ has dimension $2n-p$, then it  is a Cohen Macaulay complex space.
\begin{proof}
This is a local question at $X\times X$. Let $Z=V(\sI^2(f))$. At strict double points $(x,x')$ of $f$, the stalk $\sI^2(f)_{(x,x')}$ agrees with $\sP(f,2)_{(x,x')}$, which is generated locally by the $p$ function germs $f_j(x)-f_j(x'), 1\leq j\leq p$ (where $f_j$ is the composition of $f$ with the $j$-th coordinate function of $Y$ around $f(x)$). Thus, $Z$ is locally a complete intersection. 
Let $(x,x)\in X\times X$ be a diagonal point and denote $f_x$  the germ of $f$ at $x$. Then, $\sI^2(f)_{(x,x)}=I^2(f_x)$ and the result follows directly from Lemma \ref{lemConcciniStrickland}.\end{proof}
	\end{thm}

	\begin{thm}\label{thmI2DefinesD2}
For any finite map $f\colon X\to Y$, the ideal sheaf $\sI^2(f)$ defines the double point space $D^2(f)$.
\begin{proof}
By Proposition \ref{propUniqueMultiplePointStructure}, we only need to show that the double point structure given by $\sI^2$ satisfies conditions M1 and M2.

 To show M1, let $f$ be a stable map and denote by $Z$ the zero set of $\sI^2(f)$. By Theorem \ref{thmD2CohenMacaulay}, $Z$ is a Cohen Macaulay space of dimension $2n-p$. Now we claim that $Z$ is smooth out of the set $C=\{(x,x)\mid x\in \widehat \Sigma^2(f)\}$. By Lemma \ref{lemZerosetI2}, $Z$ consists of strict double points of $f$ and diagonal points $(x,x)$, such that $f$ is singular at $x$. 

If $(x,x)$ is a diagonal point with $x\in \Sigma^1(f)$, then the stalk $\sI^2(f)_{(x,x)}$ is the double point ideal $I^2(f_x)$ of the corank 1 map germ $f_x$ defined by $f$ at $x$ and the claim follows by Proposition \ref{propGaffneyMondAgreeCorank1}. If $(x,x')$ is a strict double point of $f$, then $\sI^2(f)$ agrees with $\sP(f,2)$ locally at $(x,x')$ by Lemma \ref{lemStalksOfMondStructure}. The claim follows since $\sP(f,2)=(f\times f)^*\sI_{\Delta(Y,2)}$ and, for every stable map $f$, the restriction of $f\times f$ to $X\times X\setminus \Delta(X,2)$ is transverse to $\Delta(Y,2)$ (Proposition \ref{propStableImpliesNormalCrossings}). 

By Proposition \ref{propDimSigmaStableMap}, the dimension of $C$ is less than or equal to $n -2(p - n + 3)<2n-p$. Hence, $Z$ is a generically smooth Cohen Macaulay space and, thus, reduced. This reduces M1 to show that $Z$ is, set theoretically, the closure of the strict double points of $f$. In other words, it suffices to show that there are no irreducible components of $Z$ consisting of points $(x,x)$ with $f$ singular at $x$. Assume that there is such a component. Then, since $Z$ is Cohen Macaulay (and hence equidimensional), the dimension of this component is $2n-p$. Therefore, the dimension of the set of singular points of $f$ is at least $2n-p$, which contradicts Proposition \ref{propDimSigmaStableMap}.

To show M2, first notice that, by Lemma \ref{lemI2AEquivalenceIsomorphism}, we can take local coordinates and assume  $X=\C^n$, $Y=\C^p$ and that the unfolding is given by 
$F(s,x)=(s,f_s(x))$, $s\in \C^r$, with $f_0=f$. Therefore it suffices to show, for any point $(x,x')$, the equality 
$$
\sI^2(f)_{(x,x')}+\langle s_i,s'_i\mid 1\leq i\leq r\rangle=\sI^2(F)_{(0,x,0,x')}+\langle s_i,s'_i\mid 1\leq i\leq r\rangle,
$$ 
where both stalks are seen as ideals in $\cO_{\C^{2(r+n)},(0,x,0,x')}$.

If $x\neq x'$, then 
\begin{align*}
\sI^2(F)_{(0,x,0,x')}&=\sP(F,2)_{(0,x,0,x')}=\langle s_i-s'_i\rangle+\langle(f_s)_j(x)-(f_s)_j(x')\rangle,\\
\sI^2(f)_{(x,x')}&=\sP(f,2)_{(x,x')}=\langle f_j(x)-f_j(x')\rangle, 
\end{align*}
so these two ideals agree modulo $\langle s_i,s'_i\rangle$.

If $x=x'$, $\sI^2(F)_{(0,x,0,x')}$ is given by the sum of $\sP(F,2)_{(0,x,0,x')}$ and the ideal generated by the $n+r$-minors of some matrix $A$, satisfying 
$$
F(s,x)-F(s',x')=A(s,x,s',x')(s-s',x-x').
$$ 
The local form of the unfolding $F$ forces $A$ to be of the form 
$$A(s,x,s',x')=\left(\begin{array}{c|c}I_r & 0 \\\hline * & \alpha_{s,s'}\end{array}\right).$$ 
Taking $s,s'=0$ we see that the submatrix $\alpha_{s,s'}$ satisfies 
$$
f(x)-f(x')=\alpha_{0,0}(x,x')(x-x').
$$ 
Therefore $\sI^2(f)_{(x,x)}$ is the sum of $\sP(f,2)_{(x,x)}$ and the ideal generated by the $n$-minors of $\alpha_{(0,0)}$, which are exactly the $n+r$-minors of $A(0,x,0,x')$. Again the equality modulo $\langle s_i,s'_i\rangle$ is immediate.
 \end{proof}
	\end{thm}

\begin{open}
Find an explicit set of generators of  the ideal defining the $k$th multiple point space $D^k(f)$, when $k\ge 3$ and $f$ has corank $\ge 2$.
\end{open}

	\section{Properties of the double point space}\label{secPropsDoublePointSpace}
Putting together Theorem \ref{thmI2DefinesD2}, Theorem \ref{thmD2CohenMacaulay}, Lemma \ref{lemZerosetI2} and Proposition \ref{lemDimensionBoundMultiplePoints}, we obtain the following:
	\begin{thm}\label{DimCorrectImpliesD2(f)CM}
For any finite map $f\colon X\to Y$
\begin{enumerate}
\item Set theoretically,  $D^2(f)$ is the union of the strict double points of $f$ and the pairs $(x,x)$ such that  $f$ is singular at $x$.
\item $D^2(f)$ has dimension $\geq 2n-p$ at every point. In particular, $D^2(f)$ is empty or $\dim(D^2(f))\geq 2n-p$.
\item  If $\dim D^2(f)=2n-p$, then $D^2(f)$ is Cohen Macaulay. 
\end{enumerate}
	\end{thm}
	
The first two items of the previous theorem can be found in \cite{Laksov77ResidualIntersectionsAndTodds} and \cite{Ronga1972La-classe-duale}.

	\begin{lem}\label{lemEdimD2AndCorankGeq2}
Let $f\colon (\C^n,0)\to (\C^p,0)$ be a finite map germ with $\corank f=k\geq 2$, then the embedding dimension of its double point space is $$\edim D^2(f)=n+k.$$

\begin{proof}
We may assume that $f$ is of the form $$(x,y)\mapsto (x,f_{n-k+1}(x,y),\dots, f_{p}(x,y)),$$ with $x=x_1,\dots,x_{n-k}, y=y_1,\dots,y_k$ and $f_j\in \km^2$, where $\km$ stands for the maximal ideal of $\cO_{2n}$.
Then, the ideal $P(f,2)+\km^2$ is generated by  $n-k$ linearly independent elements in $\km/\km^2$. Now let $\alpha$ be a matrix satisfying $f(x)-f(x')=\alpha(x-x')$. The rows corresponding to the coordinate functions $f_j, n-k+1\leq j\leq p$ have all entries in $\km$. Since there are only $n-k$ remaining rows, it follows that all the $n\times n$ minors of $\alpha$ are contained in $\km^k\subseteq \km^2$.  We obtain $\edim D^2(f)=2n-(n-k)=n+k$.  
\end{proof}
	\end{lem}

	\begin{prop}\label{D2NotRegularInCorank2}
Let $f\colon X\to Y$ be stable and $n\leq p$. Set theoretically, the singular locus of $D^2(f)$ is $$\{(x,x)\in X^2\mid x\in\widehat\Sigma^2(f)\}.$$

\begin{proof}
Let $C=\{(x,x)\in X^2\mid x\in\widehat\Sigma^2(f)\}$. In the proof of Theorem \ref{thmI2DefinesD2} we have shown that $V(\sI^2(f))$ is smooth out of $C$. Since $D^2(f)=V(\sI^2(f))$,  the singular locus of $D^2(f)$ is contained in $C$. Now let $(x_0,x_0)\in C$ and let $k=\corank f_{x_0}\geq 2$. By definition, $\sI^2(f)_{(x_0,x_0)}=I^2(f_{x_0})$, where $f_{x_0}$ is the germ of $f$ at $x_0$. From \ref{lemDimensionStrictStable} we obtain $\dim D^2(f)=2n-p$, and from Lemma \ref{lemEdimD2AndCorankGeq2}, since $2n-p\leq n<n+k$, the statement follows.
 \end{proof}
	\end{prop}

Given a finite map $f\colon X\to Y$, we denote by $p\colon D^2(f)\to X$ the map obtained by restricting the projection on the first coordinate $X^2\to X$. We define the \emph{source double point space of $f$} as $D(f)=p(D^2(f))$, with the analytic structure given by the $0$-Fitting ideal sheaf of the module $p_*\cO_{D^2(f)}$ (see \cite{MondPellikaanFittingIdeals} for details on structures given by Fitting ideals).
	\begin{cor}\label{normal} For any stable map $f\colon X\to Y$ with $n\leq p$, the projection $p\colon D^2(f)\to D(f)$ is a normalization.
\begin{proof}
 By Proposition \ref{lemDimensionBoundMultiplePoints}, $D^2(f)$ is dimensionally correct. Thus, by Theorem \ref{DimCorrectImpliesD2(f)CM},  it is a Cohen Macaulay space.  By Proposition \ref{propDimSigmaStableMap} and Proposition \ref{D2NotRegularInCorank2}, the singular locus of $D^2(f)$ is empty or has dimension $n-2(p-n+2)$.  Hence, it has codimension $\geq p-n+4\geq 4$. Thus $D^2(f)$ is a normal complex space by Serre's criterion \cite[Thm. 23.8]{MatsumuraRingTheory}. Since $f$ is stable, the strict double points are dense in $D^2(f)$ and the triple points have dimension strictly smaller. Hence $p$ is genenerically one-to-one.
\end{proof}
	\end{cor}

	\begin{cor}\label{corD2NormalD2Reduced}
If $f\colon(\C^n,0)\to(\C^p,0)$ is finitely determined and $n\leq p$, then
\begin{enumerate}
\item If $2n-p\geq 2$, then $D^2(f)$ is a normalization of $D(f)$.
\item If $2n-p\geq 1$, then $D^2(f)$ is reduced.
\end{enumerate}
\begin{proof}
By the Mather-Gaffney Criterion (see \cite{WallFiniteDeterminacyOfSmoothMapGerms}), there exists a representative of $f$ defined on a open neighbourhood of the origin $U$ (denoted also by $f$), such that $f^{-1}(0)=\{0\}$ and such that the restriction $f\vert _{U\setminus\{0\}}$ is stable. Then, we have $D^2(f)=D^2(f\vert _{U\setminus\{0\}})\cup\{(0,0)\}$ and $D^2(f\vert _{U\setminus\{0\}})$ is normal by Corollary \ref{normal}. Moreover, $D^2(f)$ has dimension $2n-p$ and is Cohen-Macaulay by Theorem \ref{thmD2CohenMacaulay}.
By Serre's criterion \cite[Thm. 23.8]{MatsumuraRingTheory},  if $\dim D^2(f)\ge1$ then it is reduced, and if $\dim D^2(f)\ge2$ then it is also normal.
\end{proof}
	\end{cor}
	
	The first item of the previous corollary can be found for $n=3, p=4$ in \cite[Prop. 4.3.1]{AltintasThesis}. The two following examples justify the need of the hypothesis in the previous corollary:
	
	\begin{ex} Take $f\colon(\C^2,0)\to(\C^3,0)$ given by
	$$
	f(x,y)=(x,y^2,y(y^2+x^2)).
	$$
	Then $D^2(f)$ is the curve in $(\C^4,0)$ defined by the equations $x=x'$, $y+y'=0$ and $x^2+y^2=0$. Since $D^2(f)$ is singular, it is not normal.
	
	\end{ex}

	\begin{ex}
Let $f\colon (\C^2,0)\to(\C^4,0)$ be given by $$(x,y)\mapsto(x^2,y^2,x^3+xy,y^3+xy).$$
A straightforward calculation yields $\dim D^2(f)=0$,
but $\edim D^2(f)=4$, by Lemma \ref{lemEdimD2AndCorankGeq2}. It follows that $D^2(f)$ is singular and, since it has dimension equal to $0$, it must be non-reduced. 
	\end{ex}
	
	\begin{prop}\label{tildeC-Mac dim n-1}
Let $f\colon X\to Y$ be a finite map  with $p=n+1$. If $f$ is generically one-to-one, then $D^2(f)$ is empty or a Cohen-Macaulay space of
dimension $n-1$.
\begin{proof}
If $f$ is finite and generically one-to-one, then the dimension of the set of strict points of $f$ is $\leq n-1$. By Lemma \ref{lemSigmaFiniteMap}, the dimension of the space of pairs $(x,x)$ with $f$ singular at $x$ is $\leq n-1$. The claim follows immediately from Theorem \ref{DimCorrectImpliesD2(f)CM}.
\end{proof}
	\end{prop}

\section {Another multiple point structure}\label{secAnotherMultiplePointStructure}

 The present section is devoted to the study of a different approach to the computation of multiple points. This alternative structure, was introduced for double points by Mond in \cite{MondSomeRemarks}, where he shows that it agrees with the usual double point structure $I^2(f)$, provided that $f$ has corank 1.  We will give some criteria for the equality of both structures of double points and show one example where they disagree. Therefore, we conclude that the new structure does not satisfy the properties M1 and M2.

Recall that the ideal sheaf $\sP(f,k)$ defines the locus of points $$(x^{(1)},\dots, x^{(k)})\in X^k,$$ such that $f(x^{(1)})=f(x^{(l)})$ for all $l\leq k$. It is clear that the zeros of $\sP(f,k)$ may contain contain some points which are not what we defined as $k$-multiple points. Indeed, for any $(x^{(1)},\dots, x^{(k-1)})$ belonging to the zeros of $\sP(f,k-1)$, the point  $(x^{(1)},\dots, x^{(k-1)},x^{(k-1)})$ belongs to $\sP(f,k)$ without imposing further conditions. For instance, the zero set of $\sP(f,k)$  contains always the small diagonal $\Delta(X,k)$. It seems a good idea to erase, taking multiplicities into account, the trivial copies of the diagonal which appear in the zeros of  $\sP(f,k)$. Locally, given two subspaces $A=V(I)$ and $B=V(J)$ of $(\C^n,0)$, we have $$\overline {A\setminus B}=V(I:J^\infty).$$ In order to erase $B$ from $A$ with multiplicity, we take the zeros of the transporter ideal $$I:J=\{h\in \cO_n\mid hJ\subseteq I\},$$ (see \cite{ASingularIntroduction} for saturation and transporter ideal).  This local definition extends to the corresponding operation between sheaves which, furthermore, preserve coherence.
\begin{definition}
 For any map $f\colon X\to Y$, we define
$$\sH^k(f)=\sP(f,k):\sI_{D(X,k)}.$$
We denote by $\tilde D^k(f)$ the complex space defined by $\sH^k(f)$. If $f$ is a map germ, then we define $H^k(f)=P(f,k):I_{D(X,k)}$ and $\tilde D^k(f)=V(H^k(f))$.
\end{definition}

	\begin{lem}\label{lemVI2=VH2SetTheoretically} \cite{MondSomeRemarks} For any map germ $f$ of corank $k$:	
 \begin{enumerate}
\item$I^2(f)\subseteq H^2(f),$
\item$H^2(f)^k\subseteq I^2(f).$
\end{enumerate}
In particular, if $\corank f=1$, then $I^2(f)=H^2(f)$.  
	\end{lem}
	\begin{cor}\label{corI2SubsheafH2}
For any $f\colon X\to Y$, the spaces $D^2(f)$ and $\tilde D^2(f)$ agree set-theoretically. At the level of schemes, the ideal sheaf $\sI^2(f)$ is a subsheaf of $\sH^2(f)$ and they both agree outside the space $$\{(x,x)\in \Delta(X,2)\mid x\in \widehat\Sigma^2(f)\}.$$ 
\begin{proof}
Out of the diagonal, both $\sH^2(f)$ and $\sI^2(f)$ agree with $\sP(f,2)$. Let $(x,x)\in D(X,2)$ be a diagonal point and $f_x$ be the germ of $f$ at $x$.  We have the equalities $\sI^2(f)_{(x,x)}=I^2(f_x)$ and $\sH^2(f)_{(x,x)}=H^2(f_x)$. The result follows directly from the previous lemma.
\end{proof}
	\end{cor}
The following lemma is standard material, so we omit its proof.
	\begin{lem}\label{lemAssociatedPrimesQuotient}
Let $I,J$ ideals in a noetherian ring $R$ and let be $I=\bigcap_{i=1}^s \kq_i$
a minimal primary decomposition, so that the associated primes of $R/I$ are
$Ass(R/I)=\{\sqrt \kq_1,\dots,\sqrt \kq_s\}.$
Then, the associated primes of $R/(I:J)$ satisfy:
$$Ass(R/(I:J))\subseteq \{\sqrt{q_i}\mid J\nsubseteq \kq_i\}.$$ 
	\end{lem}
	
	\begin{ex}
The inclusion is strict in general. Let $\kq_1=\langle x^2,y^2\rangle$ and $\kq_2=\langle x^2,xy^2,y^3,z\rangle$ be primary ideals and $J=\langle x,y\rangle$. The associated primes of $I=\kq_1\cap\kq_2$ are $J=\sqrt{\kq_1}$ and $\langle x,y,z\rangle=\sqrt{\kq_2}$, since $\kq_1\not\subset\kq_2$ and $\kq_2\not\subset\kq_1$. Moreover, we have $J\not\subset\kq_2$. However, the only associated prime of $I:J$ is $J$. 
	\end{ex}

	\begin{thm}\label{ThmI2(f)=H2(f)}
If $f\colon X\to Y$ satisfies
\begin{enumerate}
\item $\dim D^2(f)=2n-p$,
\item $\dim \widehat\Sigma^2(f)<2n-p$,
\end{enumerate}
then $\sI^2(f)=\sH^2(f)$.
\begin{proof}
 By Corollary \ref{corI2SubsheafH2}, the stalks of $\sI^2(f)$ and $\sH^2(f)$ agree out of the space $\{(x,x)\in U\times U\mid x\in \widehat\Sigma^2(f)\}$, which is a space of dimension $<2n-p$.

Again by Corollary \ref{corI2SubsheafH2}, we have $\sI^2(f)\subseteq \sH^2(f)$. Therefore, the space where $\sI^2(f)$ and $\sH^2(f)$ disagree is
$$Z=\{(x,x')\in U\times U\mid \sI^2(f)_{x,x'}\subsetneq \sH^2(f)_{x,x'}\}=\supp(\sH^2(f)/\sI^2(f)).$$
Since the sheaf $\sH^2(f)/\sI^2(f)$ is coherent , we have $$Z=V(\sA nn(\sH^2(f)/\sI^2(f))).$$ The stalks of this sheaf are 
\begin{align*}
\sA nn(\sH^2(f)/\sI^2(f))_{x,x'}	&=\Ann(\sH^2(f)_{x,x'}/\sI^2(f)_{x,x'})\\
						&=\sI^2(f)_{x,x'}:\sH^2(f)_{x,x'}.
\end{align*}
 Since the zero set of a sheaf depends only of its stalks, we have the equality $Z=V(\sI^2(f):\sH^2(f)).$
We already know that $Z$ is contained in the diagonal. 
By hypothesis, $D^2(f)$ has dimension $2n-p$ and thus, by Theorem \ref{thmD2CohenMacaulay} and Theorem \ref{thmI2DefinesD2}, $D^2(f)$ is Cohen-Macaulay, and hence equidimensional. The germ of $Z$ at any point $(x,x)\in \Delta(X,2)$ is $V(\sI^2(f)_{x,x}:\sH^2(f)_{x,x})=V(I^2(f_x):H^2(f_x))$ and, by Lemma \ref{lemAssociatedPrimesQuotient}, all the associated primes of $I^2(f_x):H^2(f_x)$ are associated primes of $I^2(f_x)$.  It follows that $Z$ is equidimensional of dimension $2n-p$ or empty. Since we have shown before that $Z$ is a subspace of a space of dimension $<2n-p$, we conclude that $Z$ is empty.
\end{proof}
	\end{thm}

	\begin{cor}
Assume $p=n$.  Then $\sI^2(f)=\sH^2(f)$ for any finite map $f\colon X\to Y$ .
\begin{proof}If $D^2(f)$ is empty, then the result is trivial. Otherwise, let $n=\dim X$. Since $f$ is finite, then $\dim D^2(f)=n$ and, by Lemma \ref{lemBoundDimHatSigmaFiniteMap}, we have $\dim \widehat \Sigma^2(f)\leq n-2$.
\end{proof}
	\end{cor}

	\begin{cor}\label{corH2I2AgreeNN+1GenOneToOne}
Assume $p=n+1$. If $f\colon X\to Y$is finite and generically one-to-one, then $I^2(f)=H^2(f)$.
\begin{proof}As in the previous corollary, if $D^2(f)$ is empty, then the result is trivial. Otherwise, let $n=\dim X$. If the map $f$ is finite and generically one-to-one, then  $\dim D^2(f)=n-1$. By Lemma \ref{lemBoundDimHatSigmaFiniteMap} we have $\dim \widehat\Sigma^2(f)\leq n-2$.
\end{proof}
	\end{cor}

	\begin{cor}\label{corI2EqualsH2ForStable}
Any finite stable map $f$ satisfies $\sI^2(f)=\sH^2(f)$.
\begin{proof} If $D^2(f)=\emptyset$, then the statement follows trivially from Lemma \ref{lemVI2=VH2SetTheoretically}. If $D^2(f)\neq \emptyset$, then the stability of $f$ implies $\dim D^2(f)=2n-p$ (Lemma \ref{lemDimensionStrictStable}) and we have $\dim \widehat\Sigma^2(f)=n-2(p-n+3)<2n-p$ (Proposition \ref{propDimSigmaStableMap}).
\end{proof}

	\end{cor}

	\begin{cor}
If $f\colon(\C^n,0)\to(\C^p,0)$ is a finitely determined germ and $p<2n$, then $I^2(f)=H^2(f)$ and both ideals are reduced.
\begin{proof}
By Mather-Gaffney criterion, we can find a representative $f$ which is stable out of $\{0\}$. Then, $D^2(f)$ is reduced and, by the previous corollary, $\sI^2(f)$ and $\sH^2(f)$ agree out of $\{0\}$. Moreover, $\cO_{2n}/I^2(f)$ is Cohen Macaulay and, thus, equidimensional of dimension $2n-p>0$. As we saw in the proof of Theorem \ref{ThmI2(f)=H2(f)}, these sheaves can only differ on some zeros of associated primes of $\sI^2(f)$, which are spaces of dimension $>0$. Therefore $\sI^2(f)=\sH^2(f)$ 
\end{proof}
	\end{cor}

\begin{ex}\label{exDoubleCone}
Let be $f\colon(\C^2,0)\to(\C^3,0)$ the `Double Cone', given by  $$(x,y)\mapsto(x^2,y^2,xy).$$
A straightforward computation with {\sc Singular} yields $I^2(f)=A\cap B_1$ and $H^2(f)=A\cap B_2$, where
\begin{itemize}
\item[]$A=\langle x+x', y+y'\rangle,$
\item[]$B_1=\langle x^2,xx',xy,x'^2,x'y',y^2,yy',y'^2,xy'+x'y\rangle,$
\item[] $B_2=\langle x^2,xy,y^2,x',y'\rangle.$
\end{itemize}
The ideal $A$ defines a reduced plane, while $V(B_i)=\{0\}$. Therefore, $\dim D^2(f)=2$ and Theorem \ref{ThmI2(f)=H2(f)} does not apply here. Indeed, $D^2(f)$ has an embedded component, namely $V(B_i)$, and thus it is is not a Cohen-Macaulay space. 

Now we show  that in this situation $H^2(f)$ may not behave well under deformations:
Take the unfolding $F\colon (\C^3,0)\to(\C^4,0)$ given by
$$(t,x,y)\mapsto(t,f_t(x,y)),\ f_t(x,y)=(x^2,y^2,xy+t(y^3+x^3)).$$
For sufficiently small $t\neq 0$, $f_t$ is $\cA$-equivalent to the map in Example \ref{exDoubleFoldTrifold}. Since $f_t$ is generically one-to-one, for all $t\neq 0$ , we conclude that the map $F$ is generically one-to-one. From Corollary \ref{corH2I2AgreeNN+1GenOneToOne} it follows $I^2(F)=H^2(F)$ and, since $I^2(f)$ behaves well under deformations, we have $H^2(F)+\langle t\rangle=I^2(f)\neq H^2(f).$
	\end{ex}

\bibliography{MyBibliography} 
\bibliographystyle{plain}	

\end{document}